\documentclass[10pt,9pt]{article}
\usepackage{amssymb}
\usepackage{makeidx}
\usepackage{amsmath}
\usepackage{mathrsfs}
\usepackage{color}
\usepackage{epsfig}

\newtheorem{example}{Example}[section]
\newtheorem{definition}{Definition}[section]
\newtheorem{theorem}{Theorem}[section]

\begin{document}

Title: Using the Additive Property of Compactly Supported Cohomology Groups

Author: Ma\l gorzata Aneta Marciniak

Comments: LaTeX, 6 pages

Subj-class: Several Complex Variables

MCS-class: 32C36 (primary); 16E99 (secondary)

{\it Abstract: This article uses basic homological methods for evaluating examples of compactly supported cohomology groups of line bundles over projective curve.

\section{Introduction}

Compactly supported cohomology groups play a key role in complex analysis. Vanishing properties are related to solvability of $\overline{\partial}$ problem as well as analytic continuation problems (\cite{dw-1} and \cite{dw-2}). This article uses the K\"unneth formula and the additive property for finding examples of compactly supported cohomology groups. The additive property requires the definition of the inverse image of a sheaf, which is shorty explained in subsection \ref{sec:additive_prop}. Compactly supported cohomology groups of $\mathbb{C}^1$, $\mathbb{C}^*$, $\mathbb{C}^2\setminus \{(0,0)\}$ and $E_{-1}$ are evaluated in subsection \ref{subs:remove_point} and of $E_k$, the line bundles over  $\mathbb{P}^1$, in subsection \ref{sub:remove_proj}. An example of a trivial bundle over $\mathbb{P}^1$ can use he K\"unneth formula in place of the additive property, is presented in subsection \ref{sub:kunneth}.

\subsection{Definitions and Properties}

\begin{definition}{\bf (Compactly supported Dolbeault cohomology groups)}
Compactly supported Dolbeaut cohomology groups of the domain $D$ are the complex vector spaces:
\[
\mathfrak{H}_{c}^{p,q}(D)=\frac{\{\overline{\partial}\mbox{-closed forms with compact support of bidegree $(p,q)$ in $D$}\}}{\{\overline{\partial}\mbox{-exact forms with compact support of bidegree $(p,q)$ in $D$}\}}.
\]
\end{definition}
\bigskip
The following theorem shows the relationship between compactly supported cohomology groups and compactly supported Dolbeault cohomology groups.

\begin{theorem}{\bf (Dolbeault's Theorem, \cite{do})}
If $D$ is an open domain in the space of $n$ complex variables, $\mathscr{O}$ is the sheaf of germs of holomorphic functions on $D$, and $\mathfrak{H}_{c}^{p,q}(D)$ is compactly supported Dolbeault cohomology group of bidegree $(p,q$) for $D$. Then $H_{c}^q(D,\mathscr{O})=\mathfrak{H}_{c}^{0,q}(D)$.\,\rule[-.2ex]{.6ex}{1.8ex}
\end{theorem}

An alternate definition can be found in \cite{bredon}. Note that if $X$ is a compact manifold then $H_c^n(X, \mathscr{O})=H^n(X, \mathscr{O})$ fr all $n=0,1,2,\ldots$.

In particular, $H_c^0(X, \mathscr{O})=\mathfrak{H}_{c}^{0,0}(X)$ denotes (a linear space of) global holomorphic functions on $X$ with compact support. Note that
\[
H_c^0(X, \mathscr{O})=
\begin{cases}
0          &\text{if $X$ is noncompact,}\\
\mathbb{C} &\text{if $X$ is compact.}
\end{cases}
\]

In particular, $H_c^1(X, \mathscr{O})=\mathfrak{H}_{c}^{0,1}(X)$ can be seen as follows, where $\omega$ has a bidegree $(0,1)$:
\[
H_c^1(X, \mathscr{O})=\frac{\{\omega\; \mbox{with compact support such that} \;\overline{\partial}\omega=0\}}{\{\omega\; \mbox{with compact support such that} \;\overline{\partial}f=\omega \;\mbox{for some $f$ with compact support}\}}
\]

\subsection{The K\"unneth Formula}

The following exact sequences are obtained for each $n$ separately.
\begin{theorem}
If $X$ and $Y$ are locally compact Hausdorff spaces, with the sheaves $\mathscr{F}$ and $\mathscr{G}$ respectively and $\mathscr{F}*\mathscr{G}=0$, then the sequence \index{cohomology with compact support!K\"unneth formula}
\[
0\rightarrow \bigoplus_{p+q=n} H^p_c(X,\mathscr{F})\otimes H^q_c(Y,\mathscr{G})\rightarrow H^n_c(X\times Y,\mathscr{F}\otimes \mathscr{G})\rightarrow\bigoplus_{p+q=n+1} H^p_c(X,\mathscr{F})* H^q_c(Y,\mathscr{G})\rightarrow 0
\]
is exact.
\end{theorem}
Here $*$ denotes the free product. Note if one of the factors is torsion free then the free product is $0$.

\subsection{The additive property}\label{sec:additive_prop}

Compactly supported cohomology groups have all properties of a cohomology theory.
The ``additive" property, broadly used in the further part of the research, requires the notion of the inverse image of a sheaf.

\begin{definition}\textbf{(Inverse image)}
Let $f:A\rightarrow B$ be a map and let $\mathscr{G}$ be a sheaf on $B$ with canonical projection $\pi:\mathscr{G}\rightarrow B$. The inverse image sheaf $f^{*}\mathscr{G}$ is defined as
\[
f^{*}\mathscr{G}=\{(a,g)\in A\times \mathscr{G}:f(a)=\pi(g)\}.
\]
\end{definition}
In particular, if $f$ is a closed embedding, the following theorem holds.

\begin{theorem} \textbf{(\cite{iversen} III.7.6)}
Let $i:Y\rightarrow X$ be a closed embedding, then the following sequence
\[
\ldots\rightarrow H_c^q(X\setminus Y, \mathscr{F}) \rightarrow
H_c^q(X,\mathscr{F}) \rightarrow
H_c^q(Y, i^{*} \mathscr{F})\rightarrow
H_c^{q+1}(X\setminus Y, \mathscr{F})\rightarrow \ldots,
\]
is exact. \,\rule[-.2ex]{.6ex}{1.8ex}
\end{theorem}

\section{Examples for the Additive Property}

\subsection{Removing a point}\label{subs:remove_point}

Let us find compactly supported cohomology groups of $\mathbb{C}^1$ using the representation $\mathbb{P}^1=\mathbb{C}^1\cup \{\infty\}$
\begin{example}
Since $\,\mathbb{P}^1=\mathbb{C}^1\cup \{\infty\}$, then $X=\mathbb{P}^1$, $Y={\infty}$ and $X\setminus Y=\mathbb{C}^1$ implies the following exact sequence:
\begin{equation}
\begin{array}{c}
0\rightarrow H^0_c(\mathbb{C}^1,\mathscr{O})\rightarrow H^0_c(\mathbb{P}^1,\mathscr{O})\rightarrow H^0_c(\{\infty\},i^{*}\mathscr{O})\rightarrow H^1_c(\mathbb{C}^1,\mathscr{O})\rightarrow \\
\\
H^1_c(\mathbb{P}^1,\mathscr{O})\rightarrow H^1_c(\{\infty\},i^{*}\mathscr{O})\rightarrow 0.
\end{array}
\end{equation}
Since $H^0_c(\mathbb{C}^1,\mathscr{O})=0$, $H^1_c(\{\infty\},i^{*}\mathscr{O})=0$, and $H^0_c(\mathbb{P}^1,\mathscr{O})=\mathbb{C}$, the exact sequence can be written as follows:
\[
0\rightarrow \mathbb{C}\rightarrow H^0_c(\{\infty\},i^{*}\mathscr{O})\rightarrow H^1_c(\mathbb{C}^1,\mathscr{O}) \rightarrow 0.
\]
We need to find $H^0_c(\{\infty\},i^{*}\mathscr{O})$. In particular, the sheaf $i^{*}\mathscr{O}$ is simply $\mathscr{O}_{X}/\mathscr{O}_{X\setminus Y}$, since the point $Y=\{\infty\}$ is closed in $X$. The global functions at $\infty$ with coefficients in $\mathscr{O}_{\mathbb{P}^1}/\mathscr{O}_{\mathbb{C}^1}$ are those convergent series at $\infty$, which are $0$ in a neighborhood of $\infty$, so they are the germs of holomorphic functions of one variable at $\infty$. The exact sequence gives $H^1_c(\mathbb{C}^1,\mathscr{O}) = H^0_c(\{\infty\},i^{*}\mathscr{O})/\mathbb{C}$ in this sense, that two germs $f$ and $g$ represent distinct elements of the group if $f(\infty)\neq g(\infty)$. In other words:
\[
H^1_c(\mathbb{C}^1,\mathscr{O}) = \{\sum_{i< 0}a_i z^i, a_i\in\mathbb{C}\},
\]
where the series converges at infinity.
\end{example}

A similar procedure can be applied to $H^1_c(\mathbb{C}^{*},\mathscr{O})$.

\begin{example}
Note that $\mathbb{P}^1=\mathbb{C}^{*}\cup \{\infty\}\cup \{0\}$. We could equivalently use $\mathbb{C}^1=\mathbb{C}^{*}\cup \{0\}$ and relate to the previous result. The additive property provides the following exact sequence:

\[
0\rightarrow H^0_c(\mathbb{C}^{*},\mathscr{O})\rightarrow H^0_c(\mathbb{P}^1,\mathscr{O})\rightarrow H^0_c(\{0,\infty\},i^{*}\mathscr{O})\rightarrow H^1_c(\mathbb{C}^{*},\mathscr{O})\rightarrow H^1_c(\mathbb{P}^1,\mathscr{O})\rightarrow H^1_c(\{0,\infty\},i^{*}\mathscr{O})\rightarrow 0.
\]
After applying $H^0_c(\mathbb{C}^{*},\mathscr{O})=H^1_c(\mathbb{P}^1,\mathscr{O})=
H^1_c(\{0,\infty\},i^{*}\mathscr{O})=0 $ and $H^0_c(\mathbb{P}^1,\mathscr{O})=\mathbb{C}$, the exact sequence simplifies to:
\[
0\rightarrow\mathbb{C}\rightarrow H^0_c(\{0,\infty\},i^{*}\mathscr{O})\rightarrow H^1_c(\mathbb{C}^{*},\mathscr{O})\rightarrow 0.
\]
Since $H^0_c(\{0,\infty\},i^{*}\mathscr{O})/\mathbb{C}$ consist, we have
\[
H^1_c(\mathbb{C}^{*},\mathscr{O})=\{\sum_{i<0}a_i z^i, a_i\in\mathbb{C} \}\oplus \{\sum_{i> 0}b_j w^j, b_j\in\mathbb{C} \},
\]
where the series converge in a neighborhood of $\infty$ and $0$, respectively.
\end{example}

\begin{example}
The group $H^{*}_c(\mathbb{C}^2\setminus \{(0,0)\},\mathscr{O})$ can be found from the additive property. Let $X=\mathbb{C}^2$ and $Y=\{(0,0)\}$, then
\[
\begin{array}{c}
0\rightarrow H^0_c(\mathbb{C}^2\setminus \{(0,0)\},\mathscr{O})\rightarrow H^0_c(\mathbb{C}^2,\mathscr{O})\rightarrow H^0_c(\{(0,0)\},i^{*}\mathscr{O})\rightarrow H^1_c(\mathbb{C}^2\setminus \{(0,0)\},\mathscr{O})\rightarrow \\
\\
H^1_c(\mathbb{C}^2, \mathscr{O})\rightarrow H^1_c(\{(0,0)\},i^{*}\mathscr{O})\rightarrow  H^2_c(\mathbb{C}^2\setminus \{(0,0)\},\mathscr{O})\rightarrow H^2_c(\mathbb{C}^2,\mathscr{O})\rightarrow 0.
\end{array}
\]
The only nontrivial groups remain:
\[
0\rightarrow H^0_c(\{(0,0)\},i^{*}\mathscr{O})\rightarrow H^1_c(\mathbb{C}^2\setminus \{(0,0)\},\mathscr{O})\rightarrow 0
\]
thus
\[
H^1_c(\mathbb{C}^2\setminus \{(0,0)\},\mathscr{O})=H^0_c(\{(0,0)\},i^{*}\mathscr{O}) = \{\sum_{i,j> 0} a_{ij}z^iw^j, a_{ij}\in \mathbb{C}\},
\]
where the series converges in some neighborhood of $(0,0)$.
\end{example}

We will start with a basic example of a line bundle over $\mathbb{P}^1$.

\begin{example}
Using the fact that $E_{-1}=\mathbb{P}^2\setminus \{p\}$, we can evaluate compactly supported cohomology groups of $E_{-1}$. The additive property gives the exact sequence:
\[
\begin{array}{c}
0\rightarrow H^0_c(E_{-1},\mathscr{O})\rightarrow H^0_c(\mathbb{P}^2,\mathscr{O})\rightarrow H^0_c(\{p\},i^{*}\mathscr{O})\rightarrow H^1_c(E_{-1},\mathscr{O})\rightarrow H^1_c(\mathbb{P}^2,\mathscr{O})\rightarrow \\
\\
H^1_c(\{p\},i^{*}\mathscr{O})\rightarrow  H^2_c(E_{-1},\mathscr{O})\rightarrow H^2_c(\mathbb{P}^2,\mathscr{O})\rightarrow 0.
\end{array}
\]
Since $\mathbb{P}^2$ is compact $H^1_c( \mathbb{P}^2,\mathscr{O})=H^2_c( \mathbb{P}^2,\mathscr{O}) = 0$ and  $H^0_c( \mathbb{P}^2,\mathscr{O})=\mathbb{C}$. Then the sequences converts to:
\[
0\rightarrow H^0_c(\mathbb{P}^2,\mathscr{O})\rightarrow H^0_c(\{p\},i^{*}\mathscr{O})\rightarrow H^1_c(E_{-1},\mathscr{O})\rightarrow 0
\]
and
\[
0\rightarrow H^1_c(\{p\},i^{*}\mathscr{O})\rightarrow  H^2_c(E_{-1},\mathscr{O})\rightarrow 0.
\]
Since $H^1_c(\{p\},i^{*}\mathscr{O})=0$ because of dimensional reasons, we obtain that $H^2_c(E_{-1},\mathscr{O})=0$. The preceding exact sequence proves that $H^1_c(E_{-1},\mathscr{O})=H^0_c(\{p\},i^{*}\mathscr{O})/\mathbb{C}$, which in the terms of convergent series can be written as:
\[
H^1_c(E_{-1},\mathscr{O})=\{\sum_{(n,m)>(0,0)}a_{n,m}z^nw^m, a_{n,m}\in\mathbb{C}\},
\]
where the series converges near $(0,0)$ in the local coordinates.
\end{example}

\subsection{Removing a Projective Curve}\label{sub:remove_proj}

The total space of the line bundle $E_k$ with $k\in \mathbb{Z}$ consists of two coordinate patches $X_1\simeq \mathbb{C}^2$ and $X_2\simeq \mathbb{C}^2$ with coordinates $(z_1,w_1)$ and $(z_2,w_2)$ respectively, related on $X_1\cap X_2$ according to the rule $z_1=\frac{1}{z_2}$ and $w_1=z_2^k w_2$. It is a well known fact that  $H^1_c(E_k,\mathscr{O})=0$ for $k>0$, nevertheless, we will present how to obtain this result using the additive property.

In the previous section we found compactly supported cohomology groups of $E_{-1}$ using the representation $E_{-1}=\mathbb{P}^2\setminus \{p\}$, since $E_{-1}$ can be obtained by removing a point from the projective plane. This is not true for other line bundles over $\mathbb{P}^1$. We will consider the Hirzebruch surfaces $\mathscr{H}_k$ and the representation $E_{-k}=\mathscr{H}_k\setminus \mathbb{P}^1$. The Hirzebruch surface $\mathscr{H}_k$ consists of four coordinate patches $X_0, X_1, X_2, X_3, \simeq \mathbb{C}^2$ with $(z_j,w_j)\in X_j$ and the transition functions described below:
\[
z_1=\frac{1}{z_0}, \qquad w_1=z_0^{k}w_0
\]
\[
z_2={z_1}, \qquad w_2=\frac{1}{w_1}
\]
\[
z_3=\frac{1}{z_2},  \qquad  w_3=z_{2}^{-k}w_{2}
\]
\[
z_3={z_0},  \qquad  w_3=\frac{1}{w_0}.
\]
Let $Y_j\simeq \mathbb{P}^1$ and let us denote $i_j:Y_j\rightarrow\mathscr{H}_k$ with the following order:
\[
\mathbb{C}^1\times \mathbb{P}^1=\mathscr{H}_k\setminus Y_0
\]
\[
E_{-k}=\mathscr{H}_k\setminus Y_1
\]
\[
\mathbb{C}^1\times \mathbb{P}^1=\mathscr{H}_k\setminus Y_2
\]
\[
E_{k}=\mathscr{H}_k\setminus Y_3
\]
Here $Y_0$ is a a projective curve in $X_0\cup X_3$ described by equations $z_0=0$ and $z_3=0$. If $f_0(z_3,w_3)\in i_{0}^{*}\mathscr{O}$ then on $X_0\cap X_3$:
\[
f_0(z_3,w_3)=\sum_{(n,m)\geq (0,0)}a_{n,m}z^n_3w^m_3=\sum_{(n,m)\geq (0,0)}a_{n,m}z^n_0\frac{1}{w^m_0},
\]
which proves that $m=0$ so $f$ does not depend on $w_3$ and
\[
f_0(z_3,w_3)=\sum_{n\geq 0}a_{n}z^n_3.
\]
Thus $H^0_c(Y_0,i_{0}^{*}\mathscr{O})=\{\sum_{n\geq 0}a_{n}z_3^n, an\in\mathbb{C}\}$, where the series converges in a neighborhood of $z_0=0$.
The following exact sequence
\[
\begin{array}{c}
0\rightarrow H^0_c(\mathbb{P}^1\times \mathbb{C}^1,\mathscr{O})\rightarrow H^0_c(\mathscr{H}_k,\mathscr{O})\rightarrow H^0_c(Y_0,i_0^{*}\mathscr{O})\rightarrow H^1_c(\mathbb{P}^1\times \mathbb{C}^1,\mathscr{O})\rightarrow H^1_c(\mathscr{H}_k,\mathscr{O})\rightarrow \\
\\
H^1_c(Y_0,i_{0}^{*}\mathscr{O})\rightarrow  H^2_c(\mathbb{P}^1\times \mathbb{C}^1,\mathscr{O})\rightarrow H^2_c(\mathscr{H}_k,\mathscr{O})\rightarrow 0
\end{array}
\]
can be simplified to the sequences:
\[
0\rightarrow H^0_c(\mathscr{H}_k,\mathscr{O})\rightarrow H^0_c(Y_0,i_0^{*}\mathscr{O})\rightarrow H^1_c(\mathbb{P}^1\times \mathbb{C}^1,\mathscr{O})\rightarrow 0
\]
and
\[
0\rightarrow  H^2_c(\mathbb{P}^1\times \mathbb{C}^1,\mathscr{O})\rightarrow H^2_c(\mathscr{H}_k,\mathscr{O})\rightarrow 0.
\]

Then $\displaystyle{H^1_c(\mathbb{P}^1\times \mathbb{C}^1,\mathscr{O})=\{\sum_{n> 0}a_{n}z_3^n, \, a_n\in\mathbb{C}\}}$ and $H^2_c(\mathbb{P}^1\times \mathbb{C}^1,\mathscr{O})=0$.

\bigskip

Note that $Y_1$ is a projective curve in $X_0\cup X_1$ described in the local coordinates as $w_0=0$ and $w_1=0$. If $f_1(z_1,w_1)\in i_1^{*}\mathscr{O}$ then on $X_0\cap X_1$:
\[
f_1(z_1,w_1)=\sum_{(n,m)\geq (0,0)}a_{n,m}z^n_1w^m_1=\sum_{(n,m)\geq (0,0)}a_{n,m}\frac{1}{z^n_0}\left({z_0^kw_0}\right)^m=\sum_{(n,m)\geq (0,0)}a_{n,m}z_0^{\left(km-n\right)}w_0^m,
\]
which shows that $km-n\geq 0$. Thus $\displaystyle{H^0_c(Y_1,i_{1}^{*}\mathscr{O})=\{\sum_{(n,m)\geq (0,0)}a_{n,m}z^n_1w^m_1:km-n\geq 0, \, a_{n,m}\in\mathbb{C}\}}$.
The following exact sequence
\[
\begin{array}{c}
0\rightarrow H^0_c(E_{-k},\mathscr{O})\rightarrow H^0_c(\mathscr{H}_k,\mathscr{O})\rightarrow H^0_c(Y_1,i_1^{*}\mathscr{O})\rightarrow H^1_c(E_{-k},\mathscr{O})\rightarrow H^1_c(\mathscr{H}_k,\mathscr{O})\rightarrow \\
\\
H^1_c(Y_1,i_1^{*}\mathscr{O})\rightarrow  H^2_c(E_{-k},\mathscr{O})\rightarrow H^2_c(\mathscr{H}_k,\mathscr{O})\rightarrow 0
\end{array}
\]
can be simplified to the sequences:
\[
0\rightarrow H^0_c(\mathscr{H}_k,\mathscr{O})\rightarrow H^0_c(Y_1,i_1^{*}\mathscr{O})\rightarrow H^1_c(E_{-k},\mathscr{O})\rightarrow 0
\]
and
\[
0\rightarrow  H^2_c(E_{-k},\mathscr{O})\rightarrow H^2_c(\mathscr{H}_k,\mathscr{O})\rightarrow 0.
\]
Thus $\displaystyle{H^1_c(E_{-k},\mathscr{O})=\{\sum_{(n,m)> (0,0)}a_{n,m}z^n_1w^m_1:km-n\geq 0, \, a_{n,m}\in \mathbb{C}\}}$ and $H^2_c(E_{-k},\mathscr{O})=0$.

\bigskip

Computations for $Y_2$ (that is a submanifold in $X_1\cup X_2$) are similar to those for $Y_0$ and overall $H^0_c(Y_2,i_{2}^{*}\mathscr{O})=\{\sum_{n\geq 0}a_{n}z_1^n, \, _n\in\mathbb{C}\}$, where the series converges in a neighborhood of $z_1=0$.

\bigskip

Note that $Y_3$ is a projective curve in $X_2 \cup X_3$ described in local coordinates by $w_2=0$ and $w_3=0$. If $f_3(z_3,w_3)\in i_{3}^{*}\mathscr{O}$ then on $X_2 \cap X_3$:
\[
f_3(z_3,w_3)=\sum_{(n,m)\geq (0,0)}a_{n,m}z^n_3w^m_3=\sum_{(n,m)\geq (0,0)}a_{n,m}\frac{1}{z^n_2}\left({z_2^{-k}w_2}\right)^m=\sum_{(n,m)\geq (0,0)}a_{n,m}z_0^{\left(-km-n\right)}w_0^m,
\]
which shows that $-km-n\geq 0$ that is possible only for $m=n=0$. Thus $H^0_c(Y_2,i_{2}^{*}\mathscr{O})=\mathbb{C}$.

The following exact sequence
\[
\begin{array}{c}
0\rightarrow H^0_c(E_{k},\mathscr{O})\rightarrow H^0_c(\mathscr{H}_k,\mathscr{O})\rightarrow H^0_c(Y_3,i_3^{*}\mathscr{O})\rightarrow H^1_c(E_{k},\mathscr{O})\rightarrow H^1_c(\mathscr{H}_k,\mathscr{O})\rightarrow \\
\\
H^1_c(Y_3,i_3^{*}\mathscr{O})\rightarrow  H^2_c(E_{k},\mathscr{O})\rightarrow H^2_c(\mathscr{H}_k,\mathscr{O})\rightarrow 0
\end{array}
\]
can be simplified to the following:
\[
0\rightarrow H^0_c(\mathscr{H}_k,\mathscr{O})\rightarrow H^0_c(Y_3,i_3^{*}\mathscr{O})\rightarrow H^1_c(E_{k},\mathscr{O})\rightarrow 0
\]
and
\[
0\rightarrow  H^2_c(E_{k},\mathscr{O})\rightarrow H^2_c(\mathscr{H}_k,\mathscr{O})\rightarrow 0.
\]

\bigskip

Then $H^1_c(E_{k},\mathscr{O})=\mathbb{C}/\mathbb{C}=0$ and $H^2_c(E_{k},\mathscr{O})=0$.

\bigskip

\subsection{Example for the K\"unneth Formula}\label{sub:kunneth}

This section contains an example of a surface $\mathbb{P}^1\times\mathbb{C}^1$.

\begin{example}
Let us find $H^i_c(\mathbb{P}^1\times \mathbb{C}^1,\mathscr{O})$ for $i=0,1,2$ using the K\"unneth formula for products. Recall the following groups $H^0_c(\mathbb{P}^1,\mathscr{O})=\mathbb{C}$, $H^1_c(\mathbb{P}^1,\mathscr{O})=0$, $H^0_c(\mathbb{C}^1,\mathscr{O})=0$ and $H^1_c(\mathbb{C}^1,\mathscr{O})=\{\sum_{s\geq 0} a_{s}z^s, a_s\in \mathbb{C}\}$. Then the K\"unneth Formula for $\mathbb{P}^1\times \mathbb{C}^1$ gives the following for the first cohomology group of $\mathbb{P}^1\times \mathbb{C}^1$:
\[
0\rightarrow \bigoplus_{p+q=1} H^p_c(\mathbb{C}^1,\mathscr{O})\otimes H^q_c(\mathbb{P}^1,\mathscr{O})\rightarrow H^1_c(\mathbb{C}^1\times \mathbb{P}^1,\mathscr{O})\rightarrow\bigoplus_{p+q=2} H^p_c(\mathbb{C}^1,\mathscr{O})* H^q_c(\mathbb{P}^1,\mathscr{O})\rightarrow 0,
\]
which converts to:
\[
0\rightarrow H^1_c(\mathbb{C}^1,\mathscr{O})\bigoplus H^0_c(\mathbb{P}^1, \mathscr{O})\rightarrow H^1_c(\mathbb{C}^1\times \mathbb{P}^1,\mathscr{O})\rightarrow H^1_c(\mathbb{C}^1,\mathscr{O})* H^1_c(\mathbb{P}^1,\mathscr{O})\rightarrow 0,
\]
and gives $H^1_c(\mathbb{C}^1\times \mathbb{P}^1,\mathscr{O})=H^1_c(\mathbb{C}^1,\mathscr{O})=\{\sum_{s\geq 0} a_{s}z^s\;  \text{ and the series converges in a neighborhood of\;} 0\}$. Similarly for $H^2_c(\mathbb{C}^1\times \mathbb{P}^1,\mathscr{O})$:
\[
0\rightarrow H^1_c(\mathbb{C}^1,\mathscr{O})\otimes H^1_c(\mathbb{P}^1,\mathscr{O}) \rightarrow H^2_c(\mathbb{C}^1\times \mathbb{P}^1,\mathscr{O})\rightarrow 0,
\]
which gives $H^2_c(\mathbb{C}^1\times \mathbb{P}^1,\mathscr{O})= 0$.
\end{example}


\begin{thebibliography}{BBBC}

\bibitem{bonavero} L.~Bonavero and M.~Brion (Eds.),
{\it Geometry of toric varieties. Lectures from the Summer
School held in Grenoble},
June 19--July 7, 2000. S\'eminaires et Congres, 6.
Soci\'et\'e Math\'ematique de France, Paris (2002).

\bibitem{bredon} G. E.~Bredon, {\it Sheaf Theory}, 2 ed., Springer, (1997).

\bibitem{do} P.~Dolbeault, {\it Sur la cohomologie des vari\' et\' es analytiques complexes}
C. R. Acad. Sci. Paris, 236 (1953) 175 - 177.

\bibitem{dw-1} R.~Dwilewicz,
{\it Additive Riemann-Hilbert problem
in line bundles over $\mathbb{CP}^1$},
Canadian Math. Bull. {\bf 49} (2006), 72 - 81.


\bibitem{dw-2} R.~Dwilewicz,
{\it Holomorphic extensions in complex fiber bundles},
J.~Math.~Analysis and Appl. {\bf 322} (2006),


\bibitem{dw} R.~Dwilewicz,
{\it An analytic point of view at toric varieties},
Serdica Math. J. {\bf 33} (2007), 163 - 240.

\bibitem{ghu} B.~Gilligan, A.T.~Huckleberry, {\it Complex homogeneous manifolds with two ends},
Michigan Math. J. (1981).

\bibitem{fr} H.~Freudenthal, {\it \"Uber die Enden topologischer
R\"aume und Gruppen}, Math.~Z. {\bf 33} (1931), 692-713.

\bibitem{ful} W.~Fulton, {\it Introduction to Toric Varieties}, Princeton Univ.~Press, NJ, (1993).

\bibitem{iversen} B.~Iversen, {\it Cohomology of Sheaves},
 Springer-Verlag, (1986).

\bibitem{krantz} S.G.~Krantz, {\it Function
Theory of Several Complex Variables}, AMS Chelsea
Publishing, (2001).

\bibitem{marciniak} M.~A.~Marciniak, {\it Holomorphic Extensions in Toric Surfaces}, Journal of Geometric Analysis 2011, 1-23.

\end{thebibliography}
\end{document}